\documentclass[journal,twocolumn,singlespace]{IEEEtran}
\usepackage{amsbsy,amsmath,amssymb,epsfig,bbm,mathrsfs,multirow,amsthm,subcaption,caption,graphicx,epstopdf,multicol,lipsum,float}

\usepackage{adjustbox}
\usepackage{multicol}
\usepackage{enumerate}
\usepackage{lipsum}
\usepackage{tikz}
\usetikzlibrary{arrows,calc}
\tikzstyle{line} = [draw, -latex']

\newenvironment{myalign}{\small\par\nobreak\noindent\align}{\endalign}
\usepackage{enumitem}

\usepackage{algorithm}
\usepackage{algorithmic}

\setlist[itemize]{leftmargin=*}

\begin{document}

\title{Supplier Cooperation in Drone Delivery}

\author{\IEEEauthorblockN
{Suttinee Sawadsitang\IEEEauthorrefmark{1},
 Dusit Niyato\IEEEauthorrefmark{1}},
 Puay Siew Tan\IEEEauthorrefmark{2},
 Ping Wang\IEEEauthorrefmark{1}\\
\IEEEauthorblockA{
\IEEEauthorrefmark{1} School of Computer Science and Engineering, Nanyang Technological University\\
\IEEEauthorrefmark{2}Singapore Institute of Manufacturing Technology (SIMTech) A*STAR } \vspace{-5mm}	}

\maketitle\thispagestyle{empty}

\begin{abstract}Recently, unmanned aerial vehicles (UAVs), also known as drones, has emerged as an efficient and cost-effective solution for package delivery. Especially, drones are expected to incur lower cost, and achieve fast and environment friendly delivery. While most of existing drone research concentrates on surveillance applications, few works studied the drone package delivery planning problem. Even so, the previous works only focus on the drone delivery planning of a single supplier. In this paper, thus we propose the supplier cooperation in drone delivery (CoDD) framework. The framework considers jointly package assignment, supplier cooperation, and cost management. The objective of the framework is to help suppliers minimize and achieve fair share of the cost as well as reach a stable cooperation. The trade-off between using drones and outsourcing package delivery to a carrier is also investigated. The performance evaluation of the CoDD framework is conducted by using the Solomon benchmark suite and a real Singapore dataset which evidently confirms the practical findings.
\end{abstract}
\section{Introduction}
The popularity of drones has been increasing rapidly in many fields~\cite{drone_survey} including surveillance, entertainment, and wireless communications. The drone technology has been improved significantly in terms of travel reliability, efficiency, and energy/fuel consumption. Therefore, drones are regarded as a promising solution for package delivery~\cite{drone_delivery}. Many big businesses, e.g., Amazon and DHL, have used drones to deliver packages and merchandise items to their customers. However, drones may not be an optimal choice in many situations, and suppliers need to evaluate many factors for drone delivery. In particular, while a drone usually offers faster delivery, lower cost, less manpower, and more environment friendly than a ground-based vehicle, i.e., truck, the drone has certain constraints such as flying distance limit per trip and small delivery capacity. The flying distance limit and the location of a supplier's depot are therefore the important parameters affecting a serving area of the drone delivery. For example, if customers are located outside the serving area, the supplier cannot use a drone for package delivery of the customers. Similarly, drones cannot carry a heavy package. Therefore, suppliers need to seek for an alternative delivery mode, i.e., by outsourcing the delivery of some packages to a carrier. The cost of outsourcing the package delivery to a carrier is usually higher than that of the  drone delivery. This tradeoff needs to be optimized to minimize the cost and maximize the profit of the supplier. Moreover, multiple suppliers can cooperate and create a pool of suppliers. The cooperative suppliers share not only drones, but also depots and customers to serve. As such, the pool of drones can extend the drone serving area of one supplier. This can significantly reduce the delivery cost and improve their resource utilization.

According to the aforementioned scenario, suppliers have to address the associated challenges, i.e., (i) should suppliers cooperate and create a pool of drones or not, (ii) if the pool is created, how the suppliers fairly share the delivery cost incurred from the cooperation, (iii) to minimize the delivery cost, which customers should be served by which drone and which customers should be assigned to a carrier, and (iv) how many drones are needed and what is the routing path of each drone? To tackle these four challenges, we propose the supplier cooperation in drone delivery (CoDD) framework. The objective of the CoDD framework is to help the rational suppliers make the best decisions in terms of minimial cost and stability. The proposed CoDD framework has three decision-making components, i.e., supplier cooperation, cost management, and package assignment. In the supplier cooperation, we adopt the merge-and-split algorithm~\cite{merge} to decide on a coalition structure among suppliers. In the cost management component, the Shapley value~\cite{shapley} is applied to fairly distribute the incurred cost. In the package assignment, we propose a mix integer programming problem to optimize the trade-off between using drones and outsourcing package delivery to a carrier as well as to decide the routing paths of the drones. Furthermore, the performance evaluation of the CoDD framework is conducted by using two datasets, i.e., the Solomon benchmark suite and a real data from a Singapore logistics company. Some important findings are indicated from the numerical study. For example, the number of customers inside the serving area directly influences the cooperation decisions of the suppliers.


\section{Related Work}


Package delivery and vehicle routing problem has caught high attention from researchers for many years~\cite{vrp_survey}. Thanks to the recent technology innovation, a drone can be used for a variety of applications with acceptable reliability while incurring much cheaper cost than before. The authors in~\cite{UAV-refuel} studied a drone delivery problem aiming to minimize the total requirement of drone resources, i.e., fuel or battery. In this work, drones are allowed to visit multiple depots to refuel or charge their battery. The authors in~\cite{sidekick} proposed a delivery planning for a drone and a truck. Since the drone has a limited capacity, the truck serves the customers with a heavy package. They proposed two optimization problems, which are formulated as a mix integer programming. The difference of these problems is that a drone can depart and land at (i)  only the original depot and (ii) any customer location that the truck visited, where the truck acts as a mobile depot for the drone. They solved the optimization problems by a heuristic algorithms. 
The authors in~\cite{drone_delivery} considered only drones in the delivery. The factors incorporated in the planning include capacity, battery weight, changing payload weight, and reusing vehicles to reduce costs. They proposed two multi-trip drone routing problems, which have different objectives, i.e., (i) to minimize the cost subject to delivery time constraint and (ii) to minimize the time subject to the budget cost constraint. They used simulate annealing (SA) heuristic algorithm to find an optimal solution. The  authors in~\cite{k_mean} proposed the mathematical formulations for closed form estimations to address the joint truck and drone delivery problem. They used K-mean algorithm to find launch locations and adopted the genetic algorithm to solve the vehicle routing problem.  



Although there are some studies of the drone delivery, none of them considers the delivery planning for multiple suppliers. Especially, when package assignment, suppliers cooperation, and cost management are optimized jointly. Therefore, this is the focus of this paper. 

\section{System Model and Assumptions}
\begin{figure}[t]

\framebox[0.5\textwidth]{\centering
\includegraphics[width=0.5\textwidth]{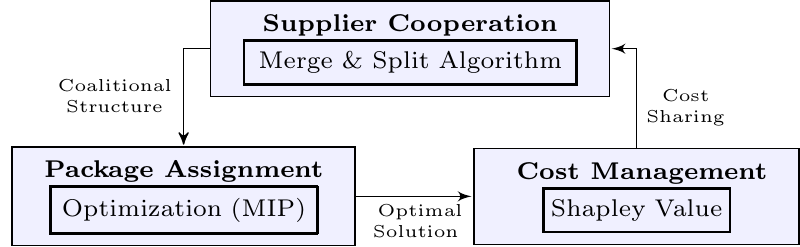}
}
\caption{The supplier cooperation in CoDD framework.}
\label{f_comp}
\vspace{-1.5em}
\end{figure}


We consider multiple suppliers using drones to deliver packages to their customers. Specifically, the suppliers have two delivery options, i.e., to use their drone or outsource the package delivery to a carrier. Therefore, the suppliers determine which packages to be delivered using the drone or carrier with different associated costs. 
The suppliers can cooperate and establish a pool of suppliers that can be shared for the  package delivery. For example, among the cooperative suppliers, a package of a customer of one supplier can be delivered by a drone belonging to the different supplier. In the following, we describe the CoDD framework, the components of which are shown in Figure~\ref{f_comp}.

\subsection{Package Assignment}
\label{sec_package}
The supplier has to decide on which drone to deliver which package, i.e., package assignment. The detail of package assignment is described as follows. Let $\mathcal{N}=\{p_1,p_2,\dots,p_{|\mathcal{N}|}\}$ be a set of all suppliers, where $|\mathcal{N}|$ is the total number of suppliers. A supplier $p \in \mathcal{N}$ has one depot. $C_p$ and $D_p$ denote the set of customers and the set of drones of supplier $p$, respectively. Two or more suppliers can cooperate and form a coalition. Let $\mathcal{P}$ be a set of suppliers in the same cooperation, i.e., coalition $\mathcal{P} \subseteq \mathcal{N}$. When suppliers cooperate, the suppliers will share their customers to serve and their drones which are referred to as a pool. For example, when suppliers $p_1$, $p_2$, and $p_3$ cooperate, the set of the customers in the pool is denoted as $\mathcal{C} = C_1 \cup C_2 \cup C_3$. Similarly, the set of the drones in the pool is denoted as $\mathcal{D} = D_1 \cup D_2 \cup D_3$. Here, $\mathcal{C}= \{c_1,c_2,\dots, c_{|\mathcal{C}|}\}$ and $\mathcal{D}= \{d_1,d_2,\dots, d_{|\mathcal{D}|}\}$ are the sets of customers and drones in the pool, where ${|\mathcal{C}|}$ and ${|\mathcal{D}|}$ are the total numbers of the customers and the drones in the pool, respectively. We use matrix $O$ to represent which customer belongs to which supplier, and the element of the matrix denoted by $o_{i,p}$ is a binary parameter. $o_{i,p} =1$ indicates that customer $i$ belongs to supplier $p$, and $o_{i,p} =0$ otherwise. Every customer in $\mathcal{C}$ has a package to be delivered. The package of customer $i$ has specific weight, which is denoted as $a_i$. A package is allowed to be transferred from one depot to another depot of cooperative supplier $p \in \mathcal{P}$ in the same pool to facilitate the drone delivery. To deliver a package of the customer, a drone must depart from the depot, which the package is located, and then come back to any depot of supplier $p \in \mathcal{P}$. Every drone $d$ has its limits, i.e., daily flying distance limit $l_d$, flying distance limit per trip $e_d$, capacity limit $f_d$, and working hour limit $h_d$. Note that the flying distance limit per trip defines the serving area. Every drone $d$ has an average flying speed, which is denoted as $s_d$. The flying distance from location $i$ to location $j$ is denoted as $k_{i,j}$. 

There are four costs involved in the package assignment including (i) the initial cost of a drone $\widehat{\mathsf{C}}_d$, such as rental fee, preparation cost, and manpower, (ii) routing cost $\bar{\mathsf{C}}_{i,j}$, which is incurred from traveling such as fuel or energy cost, (iii) the transferring cost $\ddot{\mathsf{C}}_p$, which is incurred when one supplier ships or picks up packages from other suppliers in the same pool, and (iv) outsourcing cost $\grave{\mathsf{C}}_i$ paid to a carrier.

\subsection{Supplier Cooperation}
The supplier has to decide whether or not to cooperate with other suppliers. Given the set of all suppliers $\mathcal{N}$, a certain supplier chooses to cooperate with other suppliers by forming a coalition denoted as a set $\mathcal{P}$. Let $\Phi = \{\mathcal{P}_1, \mathcal{P}_2,\dots, \mathcal{P}_{|\Phi|}\}$ denote a coalition structure, where $|\Phi|$ is the total number of coalitions in the coalition structure. The coalition structure is basically a set of all coalitions that include all the suppliers, i.e., ${\mathcal{N}} = \mathcal{P}_1 \cup \mathcal{P}_2 \cup \cdots \cup \mathcal{P}_{|\Phi|}$. To illustrate, given $\mathcal{N} = \{p_1,p_2,p_3, p_4\}$ and $\Phi = \{\{p_1,p_2\}, \{p_3,p_4\}\}$, this coalition structure consists of two coalitions, i.e., $\mathcal{P}_1 = \{p_1,p_2\}$ and $\mathcal{P}_2 = \{p_3,p_4\}$, where the supplier $p_1$ cooperates with the supplier $p_2$ and the supplier $p_3$ cooperates with the supplier $p_4$. Note that the total number of possible coalition structures can be calculated by the Bell number~\cite{shapley} based on $|\mathcal{N}|$. 


\subsection{Cost Management}
After the suppliers cooperate, it is important to achieve a fair share of cost incurred among the suppliers in the same pool, i.e., the cost management. Let $V(\mathcal{P})$ denote the delivery cost incurred from coalition $\mathcal{P}$. Let $v_i(\mathcal{P})$ denote the cost that supplier $i$ needs to pay for the delivery when joining coalition $\mathcal{P}$. Therefore, $\sum_{i\in\mathcal{P}}v_i(\mathcal{P})= V(\mathcal{P})$. We use the Shapely value as the solution of the cost management. The Shapley value~\cite{shapley} is expressed as follows:
\begin{myalign}
v_i(\mathcal{P})= \sum_{\mathcal{S} \subseteq \mathcal{P} \setminus \{i\}}\frac{|\mathcal{S}|!(|\mathcal{P}|- |\mathcal{S}| -1)!}{|\mathcal{P}|!}\left( v(\mathcal{S}\cup \{i\}) -v(\mathcal{S}) \right), \nonumber
\end{myalign}
where $|\mathcal{S}|$ denotes the total number of suppliers in set $\mathcal{S}$, and $|\mathcal{P}|$ denotes the number of suppliers in $\mathcal{P}$.

\section{Optimization Formulation}
\label{sec_formulation}
In this section, we present the optimization problem of the package assignment, which is formulated as a mix integer programming. The objective is to minimize the total delivery cost as stated in the Section~\ref{sec_package}. The objective function and the constraints are defined in (\ref{e_obj}),  and (\ref{e_con_ini}) to (\ref{e_con_tr2}), respectively. The optimization problem has six binary decision variables.
\begin{itemize}
\small
	\item $W_{d}$ indicates whether the initial cost of drone $d$ needs to be paid or not. If the payment is required, i.e., the drone will be used in the delivery, then $W_{d}=1$, and $W_{d}=0$ otherwise.
	\item $Y_{i,d,p,q}$ is an allocation binary variable of drone $d$. $Y_{i,d,p,q}=1$ if drone $d$ departs from the depot of supplier $p$ to serve customer $i$ and continues traveling to the depot of supplier $q$, and $Y_{i,d,p,q}=0$ otherwise.
	\item $Z_{i}$ indicates whether customer $i$ will be served by an outsourcing carrier or not. If $Z_{i}=1$, the package of customer $i$ is outsourced to the carrier, and $Z_{i}=0$ otherwise. 
	\item $M_{i,p,q}$ is a transferring binary variable. $M_{i,p,q} = 1$ if the package of customer $i$ needs to be transferred from the depot of supplier $p$ to the depot of supplier $q$. 
	\item $T_p$ indicates whether the transferring cost needs to be paid or not. If $T_{p} =1$, the supplier $p$ transfers its packages or picks up the packages of other suppliers in the same pool, and $T_{p} =0$ otherwise. Note that the transferring cost is not paid by for an individual supplier. It is included in the delivery cost of the pool.  
	\item $B_{p,d}$ is an auxiliary variable for eliminating impractical route. $B_{p,d}=1$ if drone $d$ departs from the depot of supplier $p$ to serve a customer and comes back to the same depot. 
\end{itemize}

\noindent Minimize:

\begin{myalign}
&\hspace{-1.5em}\sum_{d\in \mathcal{D}} \widehat{\mathsf{C}}_dW_{d} 
+ \hspace{-2em}\sum_{\substack{i \in \mathcal{C}, d\in \mathcal{D}, p, q \in \mathcal{P}}}\hspace{-2.3em}\left( \bar{\mathsf{C}}_{p,i} + \bar{\mathsf{C}}_{i,q} \right) Y_{i,d,p,q} 
+ \hspace{-0.5em}\sum_{p \in \mathcal{P}}\hspace{-0.25em}\ddot{\mathsf{C}}_{p}T_{p} + \sum_{i \in \mathcal{C}}\grave{\mathsf{C}}_iZ_i
\label{e_obj}
\end{myalign}
subject to (\ref{e_con_ini}) to (\ref{e_con_tr2}).
\begin{myalign}
&\sum_{i \in \mathcal{C}, p \in \mathcal{P}, q \in \mathcal{P} }\hspace{-1.5em}Y_{i,d,p,q} \leq \Delta W_{d}, & \forall d \in \mathcal{D}
\label{e_con_ini}
\\
&\sum_{d \in \mathcal{D}, p \in \mathcal{P}, q \in \mathcal{P} }\hspace{-1.5em}Y_{i,d,p,q} + Z_i = 1, & \forall i \in \mathcal{C}
\label{e_con_allo}
\\
&\sum_{i \in \mathcal{C}, p \in \mathcal{P}}\hspace{-1em}Y_{i,d,p,q} = \hspace{-1em}\sum_{i \in \mathcal{C}, p \in \mathcal{P}}\hspace{-1em}Y_{i,d,q,p}, & \forall q \in \mathcal{P}, d \in \mathcal{D}
\label{e_con_route}
\\
&\sum_{i \in \mathcal{C}}Y_{i,d,p,p} \leq \Delta B_{p,d}, & \forall d \in \mathcal{D}, p \in \mathcal{P}
\label{e_con_cut}
\\
& \hspace{-1em}B_{p,d}\left( \sum_{r \in \mathcal{P}}B_{r,d} -1 \right)\leq \hspace{-1.5em}\sum_{ i \in \mathcal{C}, q\in \mathcal{P}, p \neq q }\hspace{-1.5em}\Delta Y_{i,d,p,q} ,
& \forall d \in \mathcal{D}, p \in \mathcal{P}
\label{e_con_cut2}
\\
& \sum_{p \in \mathcal{P}, q \in \mathcal{P}}\hspace{-1em}a_iY_{i,d,p,q} \leq f_d, & \forall i \in \mathcal{C}, d \in \mathcal{D}
\label{e_con_cap}
\\
&\sum_{q \in \mathcal{P}}Y_{i,d,p,q}\left( k_{p,i} + k_{i,q} \right) \leq e_d, & \forall i \in \mathcal{C}, d \in \mathcal{D}, p \in \mathcal{P}
\label{e_con_pertrip}
\\
&\sum_{\i \in \mathcal{C} q \in \mathcal{P}}\hspace{-0.8em}Y_{i,d,p,q}\left( k_{p,i} + k_{i,q} \right) \leq l_d, & \forall d \in \mathcal{D}, p \in \mathcal{P}
\label{e_con_perday}
\end{myalign}

The constraint in (\ref{e_con_ini}) ensures that if packages are assigned to a drone, the initial cost of the drone must be paid, where $\Delta$~denotes a large number.  The constraint in (\ref{e_con_allo}) ensures that every package must be served either by a drone or a carrier. The constraint in (\ref{e_con_route}) ensures that each drone has the same number of times for departing and landing at any depots. The constraints in (\ref{e_con_cut}) and (\ref{e_con_cut2}) ensure that an impractical route does not exist. An example of impractical route is when a drone is assigned to serve packages from two different depots, but the drone never flies from one depot to another.

The constraints in (\ref{e_con_cap}) to (\ref{e_con_servingtime}) control flying limits of the drone delivery. The constraint in (\ref{e_con_cap}) ensures that a package, which is assigned to a drone, does not exceed the capacity limit. The constraints in (\ref{e_con_pertrip}) and (\ref{e_con_perday}) ensure that the flying distance does not exceed the limit per trip and the limit per day, respectively. The constraint in (\ref{e_con_servingtime}) ensures that the total delivery time does not exceed the limit (e.g., 8 hours). Here, traveling time and serving time are taken into account in the constraint. Note that $t_i$ is the serving time of customer $i$.

The constraint in (\ref{e_con_t_orgin}) ensures that if there is no package transferring between depots, the drone must depart from the original depot. If the package of customer $i$ is transferred to a new depot, the drone must depart from the new depot in order to serve customer $i$ as imposed by the constraint in (\ref{e_con_t_new}). The constraints in (\ref{e_con_1}) and (\ref{e_con_0}) ensure that a package is not transferred to multiple depots or the origin depot, respectively. The constraints in (\ref{e_con_tr1}) and (\ref{e_con_tr2}) ensure that the transferring cost is paid when a supplier sends or picks up packages from the pool. Note that, in the CoDD system, a drone can visit up to three depots.

\begin{myalign}
&\sum_{i,p,q}\left(\frac{k_{i,p}+k_{i,q}}{s_d} + t_i \right)Y_{i,d,p,q} \leq h_d, & \forall i \in \mathcal{C}
\label{e_con_servingtime}
\\
&o_{i,p} - \sum_{ q \in \mathcal{P}}M_{i,p,q} \leq \Delta \sum_{d \in \mathcal{D}, q \in \mathcal{P}}Y_{i,d,p,q} & \forall i \in \mathcal{C}, p \in \mathcal{P}
\label{e_con_t_orgin}
\\
&M_{i,p,q} \leq \Delta \sum_{d \in \mathcal{D}, r \in \mathcal{P}} Y_{i,d,q,r}, & \forall i \in \mathcal{C}, p, q \in \mathcal{P}
\label{e_con_t_new}
\\
&\sum_{ q \in \mathcal{P}}M_{i,p,q} \leq 1, & \forall i \in \mathcal{C}, p \in \mathcal{P}
\label{e_con_1}
\\
&M_{i,p,p} = 0, & \forall i \in \mathcal{C}, p \in \mathcal{P}
\label{e_con_0}
\\
&\sum_{i \in \mathcal{C},q \in \mathcal{P}}M_{i,p,q} \leq \Delta T_p, & \forall p \in \mathcal{P}
\label{e_con_tr1}
\\
&\sum_{i \in \mathcal{C},q \in \mathcal{P}}M_{i,q,p} \leq \Delta T_p, & \forall p \in \mathcal{P}
\label{e_con_tr2}
\end{myalign}

\vspace{-1.5em}
\section{Coalitional Game}

A coalitional game is defined by a pair $(\mathcal{N},\Pi)$, where $\Pi$ is the mapping cost for suppliers $p \in \mathcal{N}$.
 $\Pi(\mathcal{P},\Phi)$ is the collection of $\pi_q(\mathcal{P},\Phi)$ for every supplier $q \in \mathcal{N}$, where $\pi_q(\mathcal{P},\Phi)$ denotes the cost that supplier $q$ needs to pay when the coalition structure $\Phi$ is formed, i.e., $q \in \mathcal{P}, \mathcal{P} \in \Phi$. In order to select a stable coalition structure, we adopt the coalition formation algorithm in~\cite{merge}. The algorithm builds coalitions based on the preferences of the suppliers by allowing only one member to join or leave the coalition at a time, i.e., the merge-and-split algorithm. Consider two coalitions ($\mathcal{P}_1$,$\mathcal{P}_2$) and their respective coalition structure ($\Phi_1$,$\Phi_2$), where $\mathcal{P}_1 \in \mathcal{N}, \mathcal{P}_1 \in \Phi_1 $ and $\mathcal{P}_2 \in \mathcal{N}, \mathcal{P}_2 \in \Phi_2 $. Supplier $p \in \mathcal{N}$ prefers to join the coalition $\mathcal{P}_1$ over joining coalition $\mathcal{P}_2$ if $\Lambda_p(\mathcal{P}_1,\Phi_1) \leq \Lambda_p(\mathcal{P}_2,\Phi_2)$, where $\Lambda_p$ is the preference function defined as 

\begin{myalign}
\Lambda_p(\mathcal{P},\Phi) = 	\left\lbrace \begin{array}{ll}
\pi_p(\mathcal{P},\Phi), &\text{if } \pi_q(\mathcal{P},\Phi) \leq \pi_q(\mathcal{P},\Phi \setminus \{q\}) \\
& \forall q \in \mathcal{P}\setminus \{q\} \text{ and }|\mathcal{P}|\neq 1 \\
& \text{and }\mathcal{P} \notin h(p)\\
0& \text{otherwise},
\end{array} \right.
\nonumber
\end{myalign}
\noindent where $h(p)$ is a history set of the coalition structures that supplier $p$ visited, i.e., has formed before. To find the stable coalition structure, we apply the merge-and-split algorithm as shown in Algorithm~\ref{a_coalition}.


 \begin{algorithm}[]
 \scriptsize
 \caption{ for finding a stable coalition structure}
 \label{a_coalition}

 \begin{algorithmic}[1]
 
 \renewcommand{\algorithmicrequire}{\textbf{Input:}}
 \renewcommand{\algorithmicensure}{\textbf{Output:}}
 \REQUIRE set $\mathcal{N}$ and all the input parameters of package assignment optimization problem
 \ENSURE a stable coalition structure, the individual cost for each supplier, the drone delivery planning includes drone allocation and drone routing\\
 \textit{Initialisation} : All suppliers are independent (no cooperation among suppliers)
 \STATE first statement
 \WHILE { $\Phi$ has changed}
 \STATE $\Phi'$ = the list of $\Phi$ neighboring coalition structures, which can be found by using the neighborhood discovery algorithm~\cite{merge}. 
 \FOR {every $b \in \Phi'$}
 \FOR {every $p \in \mathcal{N}$}
 
 \STATE $\mathbb{C}_o$ and $\mathbb{C}_b$ are the coalitions that $p$ belongs in coalition structure $\Phi$ and $\Phi'$, respectively.
 \FOR {every $\mathcal{P} \in \Phi$}
 \STATE $V(\mathcal{P}) \leftarrow MIP(\mathcal{P})$, $\forall \mathcal{P}\in \Phi$ { \em //solve the package assignment optimization problem}
 \STATE $\Pi(\mathcal{P},\Phi) \leftarrow V(\mathcal{P})$ {\em //solve the sharing cost}
 \ENDFOR
 \STATE Do the similar loop as in line 7 to 10 for $\Phi'$.
 \IF {$\Lambda_p(\mathbb{C}_o,\Phi) \geq \Lambda_p(\mathbb{C}_b,\Phi')$}
 \STATE $\Phi = \Phi'$; {\em // the coalition structure has changed}
 \STATE \textbf{break}; {\em// break to the while loop (line 2)}
 \ENDIF \ENDFOR
 \ENDFOR
 \ENDWHILE
 \end{algorithmic}

 \end{algorithm}

\section{Experiments and Results}

We present the performance evaluation of the CoDD framework by using  two different datasets, i.e., (i) Solomon Benchmark suite (C101)~\cite{ref_solomon} and (ii) a real data from a logistic company in Singapore. We slightly synthesize the data to suit multiple depots. We implement a GAMS script for the optimization problem and solve it by CPLEX~\cite{ref_gams}. 

We experiment the CoDD framework with the following parameter setting. Four suppliers are considered. Every supplier has one drone, and all the drones are of the same type with $l_d = 150$ km, $e_d = 10$ km, $f_d=4$ kg, $h_d =8$ hrs, and $s_d = 30$ km/hr. When two suppliers cooperate, they can use either one or two drones. The numbers of customers is set equally for all suppliers, where the customers in the set $\mathcal{C}_1=\{c_1,c_5,c_9, \dots, c_{|\mathcal{N}-3|} \}$ belong to the supplier $p_1$, $\mathcal{C}_2=\{c_2,c_6,c_{10}, \dots, c_{|\mathcal{N}-2|}\}$ belong to the supplier $p_2$, $\mathcal{C}_3=\{c_3,c_7,c_{11}, \dots, c_{|\mathcal{N}-1|}\}$ belong to the supplier $p_3$, and $\mathcal{C}_4=\{c_4,c_8,c_{12}, \dots, c_{|\mathcal{N}|}\}$ belong to the supplier $p_4$. The time of loading a package for a drone is set to be $t_i = 5$ seconds. Moreover, the weight of each package is set to be $a_i=3$ kg, and thus the cost of outsourcing package to a carrier is set accordingly to be $\grave{\mathsf{C}}_i = S\$16$ based on the service by Singapost company~\cite{ref_singpost}. The cost of transferring packages is set to be $\ddot{\mathsf{C}}_p = S\$30$ for every supplier. The traveling cost is set to be $\bar{\mathsf{C}}_{i,j}= k_{i,j}\times0.105$.

\begin{figure*}
\hspace{-1em}
\begin{tabular}{cc}
\begin{tabular}{p{0.3\textwidth}p{0.3\textwidth}}
\centering \includegraphics[width=.3\textwidth, height= 5cm]{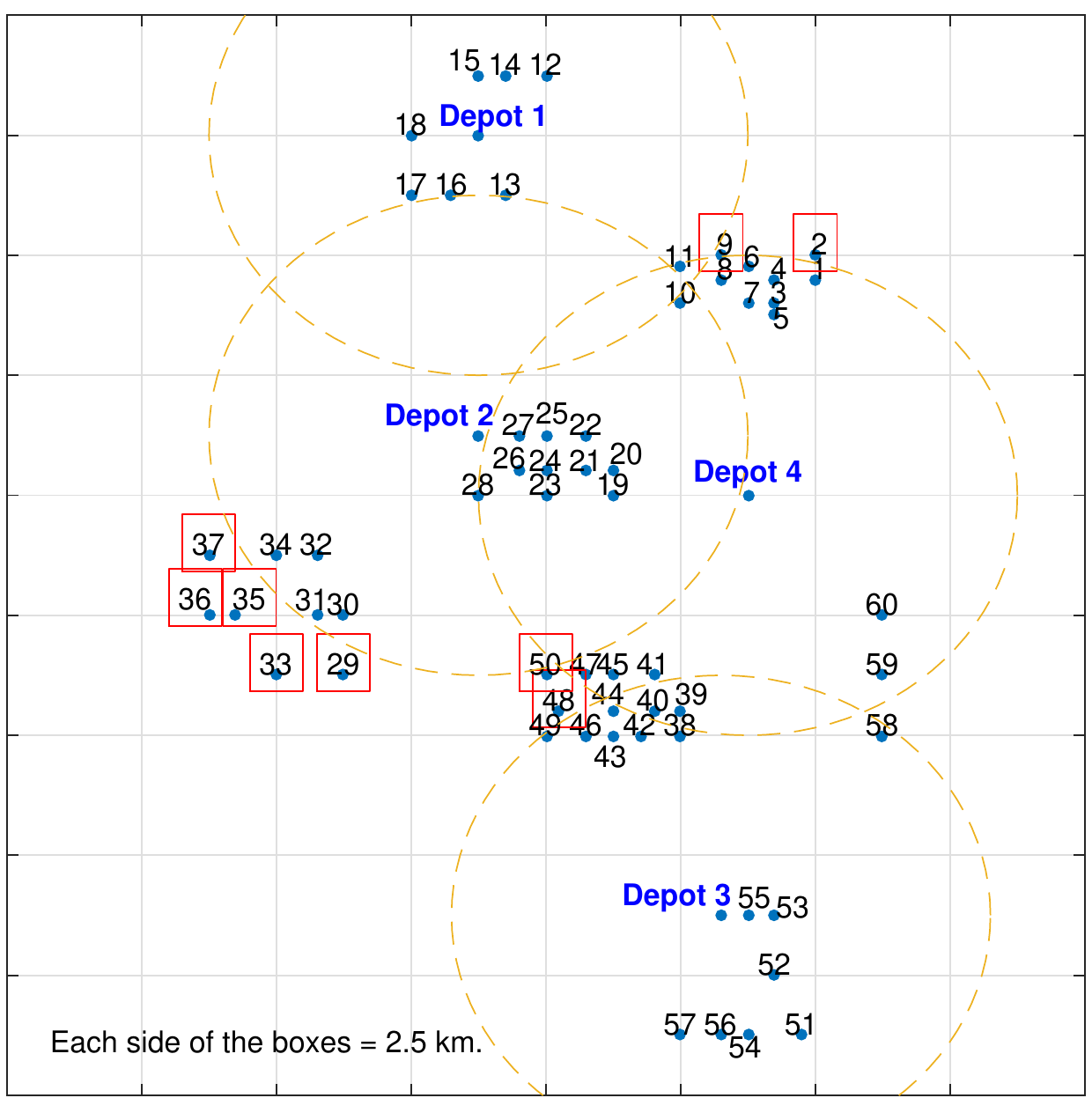}
& \includegraphics[width=.3\textwidth, height= 5cm]{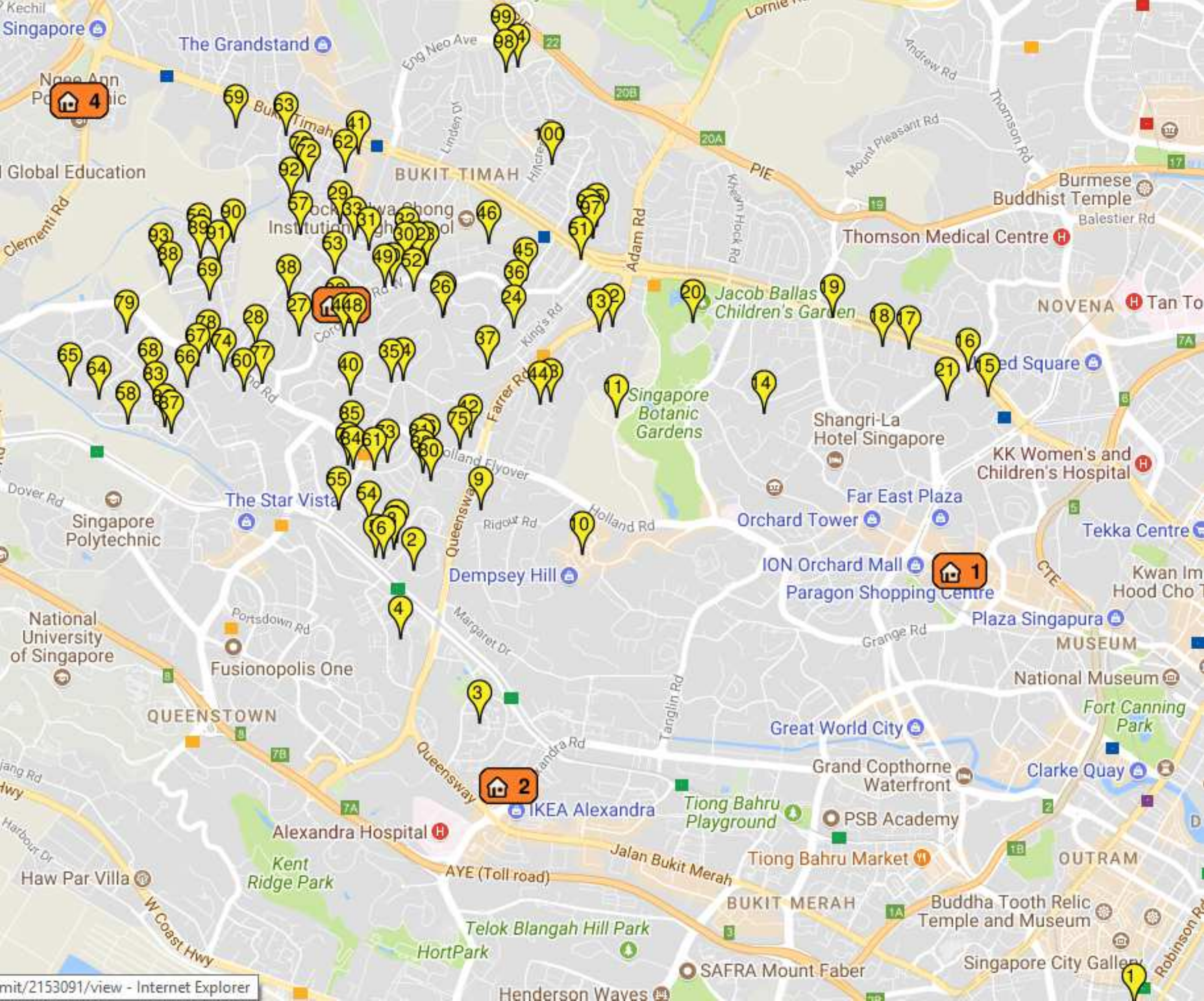}
 \\
 \scriptsize
 Depots 1, 2, 3, and 4 denote the depots of the suppliers $p_1$, $p_2$, $p_3$, and $p_4$, respectively. &
\scriptsize The orange squares represent depots. Depot3 is located very near to $c_{48}$ \\
\end{tabular} & \hspace{-2em}
\begin{tabular}{p{0.345\textwidth}}
\tiny
\hspace{-1em}\includegraphics[width=.4\textwidth]{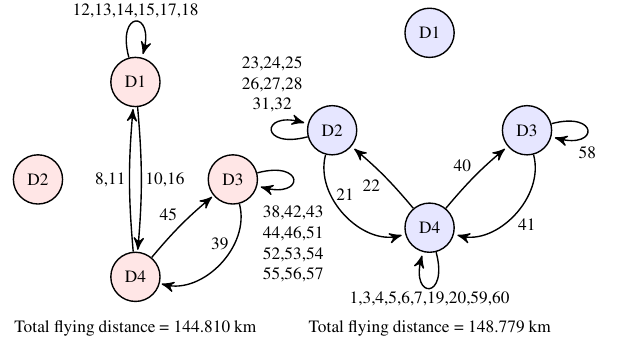}
\\
\scriptsize
The total cost is $S\$558.827$, which includes the initial cost of two drone $S\$200$, the routing cost $S\$30.827$, the outsourcing 13 packages to the carrier $S\$208$, and the transferring cost of package to the pool $S\$120$.
\end{tabular}

\end{tabular}\\
\vspace{-1em}
\begin{tabular}{p{0.6\textwidth}p{0.345\textwidth}}
\vspace{-0.8em}
\captionof{figure}{The locations of customers and depots from two data sets $(a)$-$left$ Solomon Benchmark Suite and $(b)$-$rigth$ the real Singapore data. }
\label{f_map} &  
\vspace{-0.8em}
\captionof{figure}{The result of Solomon benchmark suite.}
\label{f_solomon_shepley}
\end{tabular}
\end{figure*}

\begin{table*}
\caption{Shapley value}
\centering
\scriptsize
\begin{tabular}{|c|cccc|cccc| cccc|}
\cline{2-13}
\multicolumn{1}{c}{ }&\multicolumn{4}{|c|}{Solomon Benchmark} &\multicolumn{8}{|c|}{Real Data from a logistics company}\\ \cline{2-13}
\multicolumn{1}{c}{ }&\multicolumn{4}{|c|}{\text{Initial cost $\widehat{\mathsf{C}}_d = S\$100$ }}&\multicolumn{4}{|c|}{\text{Initial cost $\widehat{\mathsf{C}}_d = S\$100$ }} &\multicolumn{4}{|c|}{\text{Initial cost $\widehat{\mathsf{C}}_d = S\$0$ }} \\ \hline
\textbf{Coalition Structure} &\textbf{ $p_1$} & \textbf{ $p_2$} & \textbf{ $p_3$}& \textbf{ $p_4$} &\textbf{ $p_1$} & \textbf{ $p_2$} & \textbf{ $p_3$}& \textbf{ $p_4$}&\textbf{ $p_1$} & \textbf{ $p_2$} & \textbf{ $p_3$}& \textbf{ $p_4$}\\ 
\hline
$\Phi_1$= \{\{$p_1$\},\{$p_2$\},\{$p_3$\},\{$p_4$\}\} & 240.00 & 240.00 & 240.00 & 240.00 & 400.00 & 317.94 & 129.33 & 333.16 & 369.41 & 217.94 & 29.33 & 233.16 \\
$\Phi_2$= \{\{$p_1$,$p_2$\},\{$p_3$\},\{$p_4$\}\} & 205.84 & 205.84 & 240.00 & 240.00 & 367.36 & 285.31 & 129.33 & 333.16 & 302.07 & 150.60 & 29.33 & 233.16 \\
$\Phi_3$= \{\{$p_1$,$p_3$\},\{$p_2$\},\{$p_4$\}\} & 240.00 & 240.00 & 240.00 & 240.00 & 286.07 & 317.94 & 15.40 & 333.16 & 220.96 & 217.94 & -119.12 & 233.16 \\
$\Phi_4$= \{\{$p_1$,$p_4$\},\{$p_3$\},\{$p_3$\}\} & 190.51 & 240.00 & 240.00 & 190.51 & 369.06 & 317.94 & 129.33 & 302.22 & 315.81 & 217.94 & 29.33 & 179.56 \\
$\Phi_5$= \{\{$p_2$,$p_3$\},\{$p_1$\},\{$p_4$\}\} & 240.00 & 240.00 & 240.00 & 240.00 & 400.00 & 245.21 & 56.60 & 333.16 & 369.41 & 145.33 & -43.28 & 233.16 \\
$\Phi_6$= \{\{$p_2$,$p_4$\},\{$p_1$\},\{$p_3$\}\} & 240.00 & 198.50 & 240.00 & 198.50 & 400.00 & 241.87 & 129.33 & 257.09 & 369.41 & 133.80 & 29.33 & 149.02 \\
$\Phi_7$= \{\{$p_3$,$p_4$\},\{$p_1$\},\{$p_2$\}\} & 240.00 & 240.00 & 167.78 & 167.78 & 400.00 & 317.94 & 56.75 & 260.58 & 369.41 & 217.94 & -43.23 & 160.60 \\
$\Phi_8$= \{\{$p_1$,$p_2$\},\{$p_3$,$p_4$\}\} & 205.84 & 205.84 & 167.78 & 167.78 & 367.36 & 285.31 & 56.75 & 260.58 & 302.07 & 150.60 & -43.23 & 160.60 \\
$\Phi_9$= \{\{$p_1$,$p_3$\},\{$p_2$,$p_4$\}\} & 240.00 & 198.50 & 240.00 & 198.50 & 286.07 & 241.87 & 15.40 & 257.09 & 220.96 & 133.80 & -119.12 & 179.56 \\
$\Phi_{10}$= \{\{$p_1$,$p_4$\},\{$p_2$,$p_3$\}\} & 190.51 & 240.00 & 240.00 & 190.51 & 369.06 & 245.21 & 56.60 & 302.22 & 315.81 & 145.33 & -43.28 & 179.56 \\
$\Phi_{11}$= \{\{$p_1$,$p_2$,$p_3$\},\{$p_4$\}\} & 176.40 & 176.40 & 210.55 & 240.00 & 214.97 & 181.82 & 46.86 & 333.16 & 187.93 & 112.30 & -157.42 & 233.16 \\
$\Phi_{12}$= \{\{$p_1$,$p_2$,$p_4$\},\{$p_3$\}\} & 180.32 & 188.31 & 240.00 & 172.98 & 320.53 & 193.35 & 129.33 & 210.26 & 248.02 & 66.01 & 29.33 & 94.97 \\
$\Phi_{13}$= \{\{$p_1$,$p_3$,$p_4$\},\{$p_2$\}\} & 185.42 & 240.00 & 162.69 & 113.20 & 260.54 & 317.94 & -51.77 & 235.05 & \textbf{188.15} & \textbf{217.94} & \textbf{-170.89} & \textbf{127.79} \\
$\Phi_{14}$= \{\{$p_2$,$p_3$,$p_4$\},\{$p_1$\}\} & 240.00 & 192.87 & 162.16 & 120.66 & 400.00 & 206.38 & 21.25 & 221.75 & 369.41 & 107.17 & -69.86 & 122.44 \\
\textbf{$\Phi_{15}^*$ = \{\{$p_1$,$p_2$,$p_3$,$p_4$\}\}} & \textbf{156.31} & \textbf{163.79} & \textbf{138.17} & \textbf{100.57} & \textbf{231.97} & \textbf{177.81} & \textbf{-67.31} & \textbf{188.62} & 162.84 & 81.86 & -155.03 & 97.35 \\ 
\hline
\end{tabular} 
\vspace{-2em}
\label{t_shapley}
\end{table*}
\subsection{Solomon Benchmark Suite}

We evaluate the CoDD framework with 60 customers as presented in Figure~\ref{f_map}$(a)$. We set the initial cost of drones to be $\widehat{\mathsf{C}}_d=S\$100$. The solution is presented in Figure~\ref{f_solomon_shepley}, where all the suppliers cooperate ($\Phi_{15}$). Two drones are used to serve 37 customers. The routing paths are presented in Figure~\ref{f_solomon_shepley}. The solution uses only two drones instead of three drones because 9 out of 13 customers (i.e., 60-37=13) are outside the serving area, and only four customers are not worth to use one more drone. As a result, 13 customers are served by the carrier. The individual cost that each supplier needs to pay can be found in Table~\ref{t_shapley}. 


If the suppliers do not cooperate ($\Phi_1$), the suppliers need to pay $S\$240$ for outsourcing 15 packages to the carrier. The suppliers should not use a drone, since they have few customers inside their serving area, e.g., the supplier $p_1$ has only $c_{13}$ and $c_{17}$ inside the serving area.

We vary the initial cost of the drones from $\widehat{\mathsf{C}}_d = S\$100$ to $\widehat{\mathsf{C}}_d =S\$0$ to study the impact on the supplier cooperation. $\widehat{\mathsf{C}}_d =S\$0$ can represent the case when the supplier owns a drone and does not need to pay any cost for using it. However, the change from the initial cost does not affect the suppliers cooperate ($\Phi_{15}$) in this case. The reason is that the numbers of customers in the serving area of the suppliers are not large when they do not cooperate. Consequently, the cooperation between two or three suppliers still gains less benefit than that of the cooperation of all the suppliers. 

Due to the same reason, there is no change to the suppliers' cooperation when we move the depot of the supplier $p_1$ to the location which is 20 km north of the original location. The stable coalition is still  $\Phi_{15}$. 

\subsection{Real Data from a logistics company in Singapore}

We next consider 100 customers from the real Singapore data with two cases including when the initial cost is $\widehat{\mathsf{C}}_d = S\$100$ and $\widehat{\mathsf{C}}_d = S\$0$. The locations of the depots and customers are presented in Figure~\ref{f_map}$(b)$. According to the Shapley value shown in Table~\ref{t_shapley}, any supplier that cooperates with the supplier $p_3$ always achieves lower cost than that without cooperation. The reason is that the location of the depot of supplier $p_3$ can cover many customers within the drone flying range. Additionally, supplier $p_3$ achieves a negative cost in many coalition structures. The negative cost means that the supplier gains a positive revenue from the cooperation. This situation can happen, for example, when other suppliers use the drone or the depot of this supplier frequently, and thus have to share the cost by paying to the supplier. On the other hand, the location of the depot of supplier $p_1$ can cover few customers. As a result, when the initial cost is high, i.e., $\widehat{\mathsf{C}}_d = S\$100$, supplier $p_1$ does not use a drone. The supplier $p_1$ outsources all packages to the carrier. 

According to Table~\ref{t_shapley}, when the initial cost is high, i.e., $\widehat{\mathsf{C}}_d = S\$100$, the solution is that all suppliers cooperate, i.e., $\Phi_{15}$. The reason is that $\Phi_{15}$ gives the lowest cost for every supplier. When there is not initial cost, i.e., $\widehat{\mathsf{C}}_d = S\$0$, the suppliers $p_1$, $p_3$, and $p_4$ cooperate, i.e., $\Phi_{13}$, without the supplier $p_2$. This is due to the fact that suppliers $p_1$ and $p_3$ will always want to cooperate as they achieve very low cost ($\Phi_{3}$), and three suppliers cooperate will further reduce the cost comparing to the cooperation of the two suppliers. Therefore, suppliers $p_1$ and $p_3$ allow one more supplier to join in. Note that the algorithm chooses supplier $p_4$ to join before supplier $p_2$. However, after suppliers $p_1$, $p_3$, and $p_4$ are in the coalition, supplier $p_3$ does not allow supplier $p_2$ to join the coalition. Supplier $p_3$ will suffer from the higher cost if $p_2$ joins the coalition because (i) the number of customers in the serving area of pool $\Phi_{13}$ is large enough and (ii) the overlapping serving area of $p_2$ and $p_3$ decreases the resource utilization of $p_3$.
In summary, when the number of customers in the serving area is large enough, the initial cost of drones and the locations of depots can affect the stable coalition structure. 
\vspace{0em}
\section{Conclusions}
The supplier cooperation in drone delivery (CoDD) framework has been proposed to help suppliers decide whether to cooperate with other suppliers or not. 
We have proposed the optimization to address the trade-off between drone delivery and outsourcing packages to a carrier.
We have adopted the merge-and-split algorithm and Shapley value to decide how suppliers should cooperate and share the cost incurred from the delivery and cooperation. Moreover, we have performed the evaluation of the framework by using data from Solomon Benchmark suite and Singapore real data. The experiment shows that if the number of customers in the serving area is small, then outsourcing packages to a carrier is cheaper than drone delivery. When the number of customers in the serving area is large, the initial cost of drones and the locations of depots have an impact to the stable coalition structure. For the future direction, the uncertainty in drone delivery and the uncertainty in supplier behaviors will be considered. 

\section{Acknowledgment}
{\small This work was partially supported by Singapore Institute of Manufacturing Technology-Nanyang Technological University (SIMTech-NTU) Joint Laboratory and Collaborative research Programme on Complex Systems.}

\end{document}